\UseRawInputEncoding
\documentclass[twoside,12pt]{article}
\usepackage{amssymb,amsmath,bm,euscript}
\usepackage{graphicx,color,caption}
\usepackage{ifpdf}

\usepackage[colorlinks=true]{hyperref}

\usepackage{braket}

\usepackage{algorithm, algorithmic}

\newtheorem{proposition}{Proposition}[section]
\newtheorem{theorem}[proposition]{Theorem}
\newtheorem{lemma}[proposition]{Lemma}

\newtheorem{definition}[proposition]{Definition}

\newcommand{\qed}{\hphantom{.}\hfill $\Box$\medbreak}

\def\S{\mathcal{S}}
\def\P{\mathcal{P}}
\def\O{\mathcal{O}}
\def\R{\mathbb{R}}
\def\EE{\mathbb{E}}
\def\F{\mathcal{F}}
\def\I{\mathcal{I}}

\def\Q{\mathcal{Q}}
\def\A{{\mathcal{A}}}
\def\B{\mathcal{B}}
\def\C{\mathcal{C}}
\def\F{\mathcal{F}}
\def\G{\mathcal{G}}
\def\CC{\mathbb{C}}
\def\KK{\mathcal{K}}

\def\U{\mathcal{U}}
\def\V{\mathcal{V}}
\def\S{\mathcal{S}}

\def\RR{\mathcal{R}}
\def\X{{\mathcal{X}}}
\def\Y{{\mathcal{Y}}}

\def\W{\mathcal{W}}
\def\O{{\mathcal{O}}}

\def\0{{\bf 0}}

\title{\bf{T-Singular Values and T-Sketching for Third Order Tensors}}
\author{ \hspace{1mm} Liqun Qi\thanks{
Department of Applied
    Mathematics, The Hong Kong Polytechnic University, Hung Hom,
    Kowloon, Hong Kong, China; ({\tt liqun.qi@polyu.edu.hk}).}
 \ and \
 Gaohang Yu\thanks{Department of Mathematics, Hangzhou Dianzi University, Hangzhou, 310018, China; ({\tt
maghyu@hdu.edu.cn}).  This author's work was supported by National Natural Science Foundation of China (No. 12071104) and Natural Science Foundation of Zhejiang Province (No. LD19A010002).}
}

\begin{document}
\date{\today}
\maketitle

\begin{abstract}
Based upon the T-SVD (tensor SVD) of third order tensors, introduced by Kilmer and her collaborators, we
define T-singular values of third order tensors.   T-singular values of third order tensors are nonnegative scalars.   The number of nonzero T-singular values is the tensor tubal rank of the tensor.  We then use T-singular values to define the tail energy of a third order tensor, and apply it to the error estimation of a tensor sketching algorithm for low rank tensor approximation.  Numerical experiments on real world data show that our algorithm is efficient.


\end{abstract}


\section{Introduction}

Suppose that we have an input tensor $\A \in \R^{m \times n \times p}$.  Let $r$ be a target rank such that $r << \min \{ m, n \}$.   Suppose that the best rank-$r$ approximation to $\A$ is $\tilde \A$ under a certain tensor rank.   We aim to produce a low-rank approximation $\hat \A$ by sketching such that it is comparable to $\tilde A$ in the sense that
\begin{equation} \label{e1.1}
\|\A - \hat \A\|_F \approx \|\A - \tilde \A\|_F,
\end{equation}
where $\|\cdot \|_F$ is the Frobenius norm.

The tensor rank theory is still not matured yet.   Two main tensor ranks are the CP rank and the Tucker rank.    For them, the Eckart-Young theorem is unknown.   Thus, it is not easy to analyze the properties of the best rank-$r$ approximation $\tilde A$ under these two kinds of tensor ranks.  On the other hand, the Eckart-Young theorem holds for the tensor tubal rank \cite{KM11, ZEAHK14}.   See the discussion on this in Section 3 of this paper.  This motivates us to use the tensor tubal rank and T-product (tensor-tensor product) for this approach.

 Select sketch size parameters $k$ and $l$.  Draw independent standard normal tensors $\B \in \R^{n \times k \times p}$ and $\C \in \R^{l \times m \times p}$.   Then we may realize the randomised sketch $(\Y, \W)$.
\begin{equation} \label{e1.2}
\Y := \A * \B \ {\rm and}\ \W = \C * \A,
\end{equation}
where $\Y \in \R^{m \times k \times p}$, $\W \in \R^{l \times n \times p}$, $*$ is the T-product operation.    See the next section for the T-product operation.

Then we may follow the matrix approach in \cite{TYUC17} to find $\hat \A$, by the following three steps.

1. Form a T-orthogonal-triangular factorization
\begin{equation} \label{e1.3}
\Y = \Q * \RR,
\end{equation}
where $\Q \in \R^{m \times k \times p}$ is a partially orthogonal tensor and $\RR \in \R^{k \times k \times p}$ is a f-upper triangular tensor, in the sense of the T-product operation.

2. Solve a least-squares problem to find $\X £º= (\C*\Q)^\dagger * \W \in \R^{k \times n \times p}$.

3. Construct the tensor tubal rank-$k$ approximation
\begin{equation} \label{e1.4}
\hat \A := \Q*\X \in \R^{m \times n \times p}.
\end{equation}

We now need to consider several problems:

1. Is such an approach useful?

2. Is such an approach workable?

3. What is the cost of this approach?

4. What are its merits and shortcomings compared with the other approaches.

In Subsection 1.3 of \cite{TYUC17}, three situations are listed for the necessity to apply matrix sketching.   The same situations can be drawn here for tensor sketching.

To make error estimation of the proposed tensor sketching algorithm, we need to define the tail energy of a third order tensor.    For a matrix $A$, the $j$th tail energy is defined as the square root of the sum of the squares of the $i$th largest singular values of $A$ for $i > j$.  
Based upon the T-SVD (tensor SVD) factorization of third order tensors, introduced by Kilmer and her collaborators, we
define T-singular values of third order tensors.   T-singular values of third order tensors are nonnegative scalars.   The number of nonzero T-singular values is the tensor tubal rank of the tensor.  We then use T-singular values to define the tail energy of a third order tensor, and apply it to the error estimation of the tensor sketching algorithm proposed above.

In the next section, we review some preliminary knowledge on T-product and T-SVD factorization.   T-singular values of third order tensors are introduced in Section 3.   In Section 4, we propose a tensor sketching algorithm.   A probabilistic error bound is established in Section 5.   Numerical experiments are presented in Section 6.  They show that our algorithm is efficient, in particular for decay spectrum problems.

\subsection{Related Works}

{\bf T-product and T-SVD}  Kilmer and her collaborators proposed T-product and T-SVD factorization of third order tensors \cite{KBHH13, KM11, KMP08, SHKM14, ZSKA18, ZA17, ZEAHK14}.   Works on applications of T-product and T-SVD factorization include \cite{CXZ20, LYQX20, MQW20, SNZ21, XCGZ21, XCGZ21a, YHHH16, ZLLZ18}.  It is shown that they are very useful in applications.  On the other hand,
the definition of tensor tubal rank in \cite{ZEAHK14} has not explained clearly that such a definition is independent from a particular T-SVD factorization of the third order tensor.  In Definition II.7 of \cite{ZA17}, singular values were defined, and cited to \cite{Br10, KBHH13}.   We have not found the definition of singular values in \cite{Br10, KBHH13}.    Suppose that $\A \in \R^{m \times n \times p}$ have a T-SVD factorization
$$\A = \U * \S * \V^\top,$$
where $\U \in \R^{m \times m \times p}$ and $\V \in \R^{n \times n \times p}$ are orthogonal tensors,  $\S \in \R^{m \times n \times p}$ is a f-diagonal tensors, and $*$ is the t-product operation.   See the next section for the definition of orthogonal tensors, f-diagonal tensors and the t-product operation.   In Definition II.7 of \cite{ZA17}, the entries of $\S$ are called the singular values of $\A$.   This definition has several problems.

1. First, $\S$ has $mnp$ entries.   Are they all called singular values of $\A$?  These are too many.   Maybe just call the diagonal entries of the frontal slices of $\S$ singular values of $\A$.   Then there are $p\min \{ m, n \}$ singular values of $\A$.  These are still too many.   The tensor tubal rank of $\A$ is at most $\min \{ m, n \}$.  Hence, this is not consistent with the matrix case, where the number of singular values is the maximum possible rank of a matrix, and the number of nonzero singular values of a matrix is equal to its rank.

2. The values of the entries of $\S$ are dependent upon the particular T-SVD factorization.   On the other hand, the singular values of a matrix is independent from a particular SVD factorization of that matrix.   They are dependent upon the matrix only.

3. Furthermore, the diagonal entries of the frontal slices of $\S$ may be all negative.  Let $\A$ be a f-diagonal tensor such that the diagonal entries of its frontal slices are all $-1$.   Then we may let $\S = \A$, $\U$ and $\V$ be identity tensors $\I_{mmp}$ and $\I_{nnp}$ respectively.  See the next section for the definition of identity tensors.   Then this is a T-SVD factorization of $\A$.   By the modified version of Definition II.7 of \cite{ZA17} (only consider the diagonal entries of the frontal slices), all the singular values of $\A$ are $-1$.   By Definition II.16 of \cite{ZA17}, the largest singular value of $\A$ is called the tensor spectral norm of $\A$.   This would result that the tensor spectral norm of $\A$ is $-1$.

Therefore, in this paper, we define T-singular values of $\A$, such that

1. They are independent from a particular T-SVD factorization of the tensor $\A$, i.e., they are dependent to the tensor $\A$ itself.

2. There are at most $\min \{ m ,n \}$ nonzero T-singular values of $\A$.

3.  All the T-singular values are nonnegative.

4. The number of the nonzero T-singular values of $\A$ is the tensor tubal rank of $\A$.

{\bf Tensor Sketching} This is a newly developed area of tensor computation and applications \cite{ANW14, DSSW18, MB18, Pa13, PP13, WTSA15}.  There is no tensor sketching method based upon T-product operations yet.

{\bf Random Tensor Methods Based on T-product} Such a method has appeared in \cite{ZSKA18}.  The data tensor was used more than one pass.  Thus, it is a random method \cite{HMT11}, not a sketching method \cite{TYUC17}.

{\bf Matrix Sketching} Our work extends the matrix sketching methods in \cite{TYUC17} to tensors via T-product and T-SVD.   As third order tensors are much more complicated than matrices and it needs to go to and back from the Fourier domain, the extension is nontrivial.

\subsection{Notations}
In this paper, matrices are denoted by capital letters $(A,B,\ldots)$, tensors by Euler script letters $(\mathcal{ A},\mathcal{B},\ldots)$, and $\mathbb{R}$ represents a real number space, $\mathbb{C}$ represents the complex number space. For a third order tensor $\mathcal{A}\in \mathbb{R}^{m\times n \times p}$, its $(i,j,k)$-th element is represented by $a_ {ijk}$, and use the Matlab notation $\mathcal{A}(i,:,:)$, $\mathcal{A}(:,i,:)$ and $\mathcal{A}(:,: ,i)$ respectively represent the $i$-th horizontal, lateral and frontal slice of the $\mathcal{A}$. The frontal slice $\mathcal{A}(:,:,i)$ is represented by $A^{(i)}$. Define $\|\mathcal{A}\|_F:=\|\mathcal{A}(:)\|_{2}=\sqrt{\langle\mathcal{A},\mathcal{A}\rangle}=\sqrt{\Sigma_{ijk}|a_{ijk }|^2}$. $\mathcal{A}^H$ and $\mathcal{A}^\dag$ respectively represent the conjugate transpose and pseudo-inverse of $\mathcal{A}$.\par
Discrete Fourier Transformation (DFT) plays a core role in the tensor-tensor product introduced later. For a tensor $\mathcal{A}\in \mathbb{R}^{m \times n \times p}$, $\bar{\mathcal{A}}\in \mathbb{C}^{m\times n \times p}$ represents the result of DFT on $\mathcal{A}$ along the third dimension. In fact, we can use the Matlab command $\bar{\mathcal{A}}={\rm{fft}}(\mathcal{A},[\ ],3)$ to directly calculate $\bar{\mathcal{A }}$, and can use the inverse DFT to calculate $\mathcal{A}$ from $\bar{\mathcal{A}}$, that is, $\mathcal{A}={\rm{ifft}}(\bar{ \mathcal{A}},[\ ],3)$. Given a tensor $\bar{\mathcal{A}}\in \mathbb{C}^{m \times n \times p}$,
\[
\bar{A}=\rm{bdiag}(\bar{\mathcal{A}})=\begin{bmatrix}
\bar{A}^{(1)} &  &  & \\
 & \bar{A}^{(2)} &  & \\
 &  & \ddots& \\
 & &  & \bar{A}^{(p)} \\
\end{bmatrix}
\]
is a block diagonal matrix of size $mp \times np$.

\section{Tensor-Tensor Product Operations}

The T-product operation, T-SVD decomposition and tensor tubal ranks were introduced by Kilmer and her collaborators in \cite{KBHH13, KM11, KMP08, ZEAHK14}.  It is now widely used in engineering \cite{CXZ20, LYQX20, SHKM14, SNZ21, XCGZ21, XCGZ21a, YHHH16, ZSKA18, ZA17, ZLLZ18}.

For a third order tensor $\A \in \R^{m \times n \times p}$, 
as in  \cite{KBHH13,KM11}, define
$${\rm bcirc}(\A):= \left(\begin{aligned} A^{(1)}\ & A^{(p)} & A^{(p-1)} & \cdots & A^{(2)}\ \\ A^{(2)} & A^{(1)} & A^{(p)} & \cdots & A^{(3)}\\ \cdot\ \ \ & \ \cdot & \cdot\ \  & \cdots & \cdot\ \ \ \\
 \cdot\ \ \ & \ \cdot & \cdot\ \  & \cdots & \cdot\ \ \ \\
 A^{(p)} & A^{(p-1)} & A^{(p-2)} & \cdots & A^{(1)} \end{aligned}\right),$$
 and bcirc$^{-1}($bcirc$(\A)):= \A$.

Various T-product structured properties of third order tensors are based upon their block circulant matrix versions.   For a third order tensor $\A \in \R^{m \times n \times p}$, its transpose
 can be defined as
 $$\A^\top = {\rm bcirc}^{-1}[({\rm birc}(\A))^\top].$$
 This will be the same as the definition in \cite{KBHH13, KM11}.  The identity tensor $\I_{nnp}$ may also be defined as
 $$\I_{nnp} = {\rm bcirc}^{-1}(I_{np}),$$
 where $I_{np}$ is the identity matrix in $\R^{np \times np}$.  

 However, a third order tensor $\S$ in $\R^{m \times n \times p}$ is f-diagonal in the sense of \cite{KBHH13, KM11} if all of its frontal slices $S^{(1)}, \cdots, S^{(p)}$ are diagonal.   In this case, bcirc$(\S)$ may not be diagonal.

  For a third order tensor $\A \in \R^{m \times n \times p}$, it is defined \cite{KM11} that
 $${\rm unfold}(\A) := \left(\begin{aligned} A^{(1)}\\ A^{(2)}\\ \cdot\ \ \\ \cdot\ \ \\ \cdot\ \ \\ A^{(p)}\end{aligned}\right) \in \R^{mp \times n},$$
and fold$($unfold$(\A)) := \A$.   For $\A \in \R^{m \times s \times p}$ and $\B \in \R^{s \times n \times p}$, the T-product of $\A$ and $\B$ is defined as
$\A * \B :=$ fold$(${bcirc$(\A)$unfold$(\B)) \in \R^{m \times n \times p}$.   Then, we see that
\begin{equation} \label{e2.1}
\A * \B = {\rm bcirc}^{-1}({\rm bcirc}(\A){\rm bcirc}(\B)).
\end{equation}
Thus, the bcirc and bcirc$^{-1}$ operations not only form a one-to-one relationship between third order tensors and block circulant matrices, but their product operation is reserved.   By \cite{KM11}, the T-product operation (\ref{e2.1}) can be done by applying the fast Fourier transform (FFT).   The computational cost for this is $O(mnsp)$ flops.

A tensor $\A \in \R^{n \times n \times p}$ has an inverse $\A^{-1} := \B \in \R^{n \times n \times p}$ if
$$\A * \B = \B * \A = \I_{nnp}.$$
If $\Q^{-1} = \Q^\top$ for $\Q \in \R^{n \times n \times p}$, then $\Q$ is called an orthogonal tensor.   If $\Q^\top * \Q = \I_{nnp}$ for $\Q \in \R^{m \times n \times p}$, then $\Q$ is called a partially orthogonal tensor.  A tensor is called f-upper triangular or f-lower triangular, respectively,  if each frontal slice is  upper triangular or lower triangular, respectively.


\begin{lemma}\cite{KM11}\label{l2.1}
Suppose that $\mathcal{A}\in \mathbb{R}^{m\times k\times p}$ and $\mathcal{B}\in \mathbb{R}^{k\times n\times p}$ are two arbitrary tensors, let $\mathcal{F}=\mathcal{A}\ast\mathcal{B}$, then the following properties hold:
\\$(1)\|\mathcal{A}\|_F^2=\frac{1}{p}\|\bar{A}\|_F^2 = {1 \over p}\sum_{i=1}^p\|\bar{A}^{(i)}\|_F^2$;
\\$(2)\mathcal{F}=\mathcal{A}\ast\mathcal{B}$ is equivalent to $\bar{F}=\bar{A}\bar{B}$.
\end{lemma}

Let $\A, \B \in \R^{m \times n \times p}$.  Then the inner product of $\A$ and $\B$ is defined as
$$\langle \A, \B \rangle = {1 \over p}{\rm trace}(\bar A^H\bar B).$$

\begin{definition} \label{d2.2}
(Gaussian random tensor)\cite{ZSKA18} A tensor $\mathcal{G}\in \mathbb{R}^{m \times n \times p}$ is called a Gaussian random tensor, if the elements of $G^{(1)}$ satisfy the standard normal distribution (i.e.,Gaussian with mean zero and variance one), and the other frontal slices are all zeros.
\end{definition}

We leave the definition of some more important concepts, namely, T-SVD factorization and tensor tubal rank, to the next section.  In fact, we will define decay T-SVD there to specify some special T-SVD factorization which satisfies the decay property.

\section{T-Singular Values}

Let $\A \in \R^{m \times n \times p}$.   By Theorem 4.1 of \cite{KM11}, $\A$ has a T-SVD as
\begin{equation} \label{e3.5}
\A = \U * \S * \V^\top,
\end{equation}
where $\U \in \R^{m \times m \times p}$ and $\V \in \R^{n \times n \times p}$ are orthogonal, $\S$ is f-diagonal.

Note that here $\S$ is real, but not necessarily nonnegative.     Let $\A$ be a f-diagonal tensor but with negative entries.   Let $\U = \I_{mmp}$ and $\V = \I_{nnp}$.    Then $\S = \A$ and is not nonnegative.

By \cite{KM11}, we may have
\begin{equation} \label{e3.6}
 \S(1, 1, k) \ge \S(2, 2, k) \ge \cdots  \ge  \S(\min \{ m,n \}, \min \{ m,n \}, k) \ge 0,
\end{equation}
for $k = 1, \cdots, p$.

\begin{definition}  (Decay T-SVD Factorization)
A T-SVD factorization (\ref{e3.5}) is called a {\bf decay T-SVD} if (\ref{e3.6}) is satisfied.
\end{definition}

By Theorem 4.1 of \cite{KM11},
for $\A \in \R^{m \times n \times p}$, such a decay T-SVD always exists.

\begin{definition}  (T-Singular Values) Suppose that $\A \in \R^{m \times n \times p}$ with a decay T-SVD (\ref{e3.5}) such that (\ref{e3.6}) is satisfied.    The $i$th largest {\bf T-singular value} of $\A$ is defined as
\begin{equation} \label{e3.7}
\sigma_i = \sqrt{\sum_{k=1}^p \S(i, i, k)^2},
\end{equation}
for $i = 1, \cdots, \min \{ m, n \}$.
\end{definition}

In the following, we will show that the above definition of T-singular values is independent from a particular decay T-SVD form.   Before doing this, we have to discuss the tensor tubal rank of third order tensors.

In Theorem 4.3 of Kilmer and Martin \cite{KM11}, a set
\begin{equation} \label{e3.7}
M = \{ \C = \X * \Y : \X \in \R^{m \times k \times p}, \Y \in \R^{k \times n \times p} \},
\end{equation}
was defined.   This implicitly defined a rank for third order tensors in $\R^{m \times n \times p}$.
Later,  in \cite{ZEAHK14}, {\bf tensor tubal rank} was defined as the number of non-zero tubes of $\S$ in the T-SVD factorization.   It was not discussed in \cite{ZEAHK14} that the definition of the tensor tubal rank is independent from a particular T-SVD factorization of the third-order tensor.  In \cite{ZLLZ18}, it was proved that for $\A \in \R^{m \times n \times p}$, its tensor tubal rank $r$ is equal to the smallest integer $r$ such that
$$\A = \X * \Y,$$
where $\X \in \R^{m \times r \times p}$ and $\Y \in \R^{r \times n \times p}$.   This is independent from a particular T-SVD factorization of $\A$.   Thus, we may formally have the definition of tensor tubal rank as follows.

\begin{definition}
(Tensor $\rm{tubal}$ rank): The $\rm{tubal}$ rank $\rm{rank}_t(\mathcal{A})$ of tensor $\mathcal{A}\in \mathbb{R}^{m\times n\times p}$ is defined as the smallest integer $r$ such that
$$\A = \X * \Y,$$
where $\X \in \R^{m \times r \times p}$ and $\Y \in \R^{r \times n \times p}$.
\end{definition}

Hence, the set $M$ defined by (\ref{e3.7}) is in fact the set of third order tensors in $\R^{m \times n \times p}$, whose tensor tubal rank is not greater than $r$.   In fact, Theorem 4.3 of \cite{KM11} is the Eckart-Young theorem for third order tensors under the tensor tubal rank.   In \cite{ZSKA18}, Theorem 4.3 was referred as a theorem for the best ``multirank-$k$'' approximation.   This is not precise.   Multirank was introduced in \cite{KBHH13} as a vector rank.  It is not clear what the best rank-$k$ approximation about a vector rank means.   However, this minor inexactness does not conceal the significance of Theorem 4.3 of \cite{KM11} as well as the other works of Kilmer and her collaborators on T-SVD factorization.  From our point of view, Theorem 4.3 of \cite{KM11} is celebrated, as it is the only version of the Eckart-Young theorem for third order tensors known until now.  What we wish to do is to define T-singular values of third order tensors to make the T-SVD and tensor tubal rank theory complete.

We have the following theorem.

\begin{theorem}   For $\A \in \R^{m \times n \times p}$, the $i$th largest T-singular value of $\A$,
defined by (\ref{e3.7}), is independent from the particular decay T-SVD of $\A$, i.e., it is only dependent upon $\A$ itself.   Furthermore, the number of nonzero T-singular values of $\A$ is the tensor tubal rank of $\A$.
\end{theorem}
{\bf Proof}   With the above discussion, Theorem 4.3 of \cite{KM11} actually shows that the best tensor tubal rank-$j$ approximation of $\A$ is
$$\A_j = \sum_{i=1}^j \U(:,i,:) * \S(i,i,:) * \V(:,i,:)^\top,$$
for a decay T-SVD factorization of $\A$.
By (\ref{e3.5}), we have a decay T-SVD factorization of $\A_j$ as
$$\A_j = \U * \S_j * \V^\top,$$
where $\S_j$ is f-diagonal, $\S_j(i,i,:) := \S(i,i,:)$ for $i \le j$ and $\S_j(i,i,:) := \0$ for $i > j$.
By Lemma 3.19 of \cite{KM11}, orthogonal tensors preserve the Frobenius norm.   Then (\ref{e3.7}) informs us that
$$\sigma_j^2 = \|\A_j - \A_{j-1}\|_F^2,$$
where $\A_0 \equiv \A$.    By Theorem 4.3 of \cite{KM11} and above discussion, the Frobenius norms of $\A_j$ and $\A_{j-1}$ as well as $\A_j - \A_{j-1}$ are independent from a particular T-SVD factorization of $\A$.   This also shows that the number of nonzero T-singular values of $\A$ is the tensor tubal rank of $\A$.
\qed

The relation between nonzero T-singular values and the tensor tubal rank is the same as the relation between nonzero singular values and the rank in the matrix case.   This is a very good property.

For a matrix $A \in \R^{m \times n}$ or $\CC^{m \times n}$, denote the $i$th largest singular value of $A$ by $\sigma_i(A)$.   Then the $j$th tail energy \cite{TYUC17} is defined by
\begin{equation} \label{e3.10}
\tau_j^2(A) := \min_{{\rm rank}(B) < j} \|A - B\|_F^2 = \sum_{i\ge j} \sigma_j^2(A).
\end{equation}
The equality follows from the Eckart-Young theorem for matrices.

\begin{definition}  (Tail Energy) Suppose that $\A \in \R^{m \times n \times p}$, and its $i$th largest T-singular value is $\sigma_i$ for $i = 1, \cdots, \min \{ m, n \}$.  Then we the {\bf $j$th tail energy} of $\A$ is defined as
\begin{equation} \label{e3.11}
\tau_j^2(\A) := \sum_{i\ge j} \sigma_i^2(A).
\end{equation}
\end{definition}

We have
$$\tau_j^2(\A) = \min_{{\rm rank_t}(\B) < j} \|\A - \B\|_F^2.$$
The equality follows from Theorem 4.3 of \cite{KM11}.

We have the following proposition for the tail energy of a third order tensor.   Its proof follows (3.1), (4.2) and the proof of Theorem 4.3 of \cite{KM11}.

\begin{proposition} \label{p3.4}
Suppose that $\A \in \R^{m \times n \times p}$ and $\bar \A$ is the result of DFT on $\A$ along the third dimension.   Let $j$ be a positive integer satisfying $j \le \min \{ m, n \}$.  Then
$$\tau_j^2(\A) = {1 \over p}\sum_{i=1}^p \tau_j^2\left( \bar A^{(i)}\right).$$
\end{proposition}

We define T-singular values of $\A$ such that we may define the $j$th tail energy of $\A$ by (\ref{e3.11}), for our error estimation of the proposed T-sketching algorithm.    The base of our definition is from \cite{KM11}.   We think a clear definition of T-singular values of $\A$, which explicitly extends the definition of singular values of matrices, is meaningful and useful.

\section{The T-Sketching Algorithms}


\begin{algorithm}[H]
\renewcommand{\thealgorithm}{} 
\caption{Framework of T-Sketching for Low-Rank Approximation, spatial domain version}
\label{alg:A}
\begin{algorithmic}
\STATE  {\textbf{Input:} tensor $\mathcal{A} \in \mathbb{R}^{m\times n\times p}$; sketch size parameters $k\le l$}
\STATE  {\textbf{Output:} the rank-$k$ approximation $ \hat \A $}
\STATE  (1) Generate Gaussian random tensors $\B \in \R^{n\times k\times p}$, $\mathcal{\C} \in \mathbb{R}^{l\times m\times p}$;
\STATE  (2) Form the range sketch of $\A$ as  $\Y = \A * \B$ and the co-range sketch $\W = \C * \A$;
\STATE  (3) Construct the tensor $\Q$ by using T-QR factorization for $\Y$;
\STATE  (4) Solve a least-squares problem to find $\X=(\C*\Q)^\dagger * \W \in \R^{k \times n \times p}$;
\STATE  (5) Construct the rank-$k$ approximation $ \hat \A := \Q*\X \in \R^{m \times n \times p}$.
\end{algorithmic}
\end{algorithm}


\subsection{Sketching The Input Tensor}

The storage cost for the sketch cost for the sketch $(\Y, \W)$ is $(mk+ln)p$ floating point numbers.   The storage cost for two standard normal test tensors is $(nk+lm)p$ floating-point numbers.

By Section 2, the computational cost for forming the sketch (\ref{e1.2}) is $O(mn(k+l)p)$ flops.

\begin{algorithm}[H]
\renewcommand{\thealgorithm}{1} 
\caption{T-Sketching the input tensor}
\label{alg:T-Sketch}
\begin{algorithmic}
\STATE {\textbf{Require:} Input tensor $\mathcal{A} \in \mathbb{R}^{m\times n\times p}$; sketch size parameters $k\le l$}
\STATE {\textbf{Ensure:} Constructs test tensors $\mathcal{\B} \in \mathbb{R}^{n\times k\times p}$ and $\mathcal{\C} \in \mathbb{R}^{l\times m\times p}$, range sketch $\Y = \A * \B \in \mathbb{R}^{m\times k \times p}$ and co-range sketch $\W = \C * \A\in \mathbb{R}^{l\times n\times p}$ as private variables}
\STATE {\textbf{private:} $\B,\C,\Y,\W$, \; \; \; \  \ \ \ \ \; \; \; \  \ \ \ \ \  \ \ \ \ $\rhd$ Internal variables for Sketch object}
\STATE {\textbf{function} TSKETCH$(\A;k;l)$}
\STATE  1 $\B\leftarrow\mbox{zeros}(n,k,p)$; $\B(:,:,1)\leftarrow\mbox{randn}(n,k)$;
\STATE  2 $\C\leftarrow\mbox{zeros}(l,m,p)$; $\C(:,:,1)\leftarrow\mbox{randn}(l,m)$;
\STATE  3 $\B(:,:,1)\leftarrow\mbox{orth}(\B(:,:,1))$; $\C(:,:,1)^*\leftarrow\mbox{orth}(\C(:,:,1)^*)$;\;   $\rhd$(optional) Improve numerical stability
\STATE 4 $\bar{\A}\leftarrow{\rm{fft}}(\mathcal{A},[\ ],3)$, $\bar{\B}\leftarrow\rm{fft}(\mathcal{B},[\ ],3)$, and $\bar{\C}\leftarrow\rm{fft}(\mathcal{C},[\ ],3)$
\STATE 5 \textbf{for} $i\leftarrow 1 $ \textbf{to} $p$  \textbf{do}
\STATE 6 \ \   $\bar{\Y}^{(i)} = \bar{\A}^{(i)} \bar{\B}^{(i)}$;  $\bar{\W}^{(i)} = \bar{\C}^{(i)} \bar{\A}^{(i)}$;
\STATE 7 \textbf{end}
\STATE 8  Construct the sketch $ \Y \leftarrow \rm{ifft}(\bar{\Y},[\ ],3)$; $ \W \leftarrow \rm{ifft}(\bar{\W},[\ ],3)$
\end{algorithmic}
\end{algorithm}


\subsection{The T-QR Process}

Consider the T-QR process (\ref{e1.3}) $\Y =: \Q * \RR$, where $\Y \in \R^{m \times k \times p}$, $\Q \in \R^{m \times k \times p}$ is partially orthogonal, and $\RR \in \R^{k \times k \times p}$ is f-upper triangular.   
To ask $\Q$ to be partially orthogonal, we actually demand $\Q^\top * \Q = \I_{kkp}$.   We may use the T-QR process described in \cite{KM11, ZSKA18}.   Then, what is the computational cost of this step?

\subsection{Solve A Least-Squares Problem}

Recently, in \cite{MQW20}, Miao, Qi, and Wei introduce the definition of Moore-Penrose inverse for a given tensor $\A$ such that
\begin{equation}\label{eq:Moore-Penrose-Inverse}
\A^\dagger={\rm bcirc}^{-1}(({\rm bcirc}(\A))^\dagger).
\end{equation}
Further, they point out that for the following least squares problem based on T-product,
\begin{equation}\label{eq:least-square}
\min_{\X} \frac{1}{2}\|\A*\X-\B\|_F^2,
\end{equation}
its least squares solution is $\X=\A^\dagger*\B$. According to Lemma \ref{l2.1}, the objective function in (\ref{eq:least-square}) can be reformulated as
\begin{equation}\label{eq:least-square2}
\min_{\X} \frac{1}{2p}\|\bar{A}\bar{X}-\bar{B}\|_F^2,
\end{equation}
then, $\bar{X}=\bar{A}^\dagger \bar{B}$.

\subsection{The T-Sketching Algorithms}
For the convenience of error analysis, we present an implementation of T-Sketching Algorithm in the Fourier domain.

\begin{algorithm}[H]
\renewcommand{\thealgorithm}{2} 
\caption{T-Sketching for Low-Rank Approximation, Fourier domain version}
\label{alg:B}
\begin{algorithmic}
\STATE {\textbf{Input:} $\mathcal{A} \in \mathbb{R}^{m\times n\times p}$; sketch size parameters $k\le l$}
\STATE {\textbf{Output:} low-rank approximation $ \hat \A$}

\STATE 1 Generate Gaussian random tensors $\B \in \R^{n\times k\times p}$, $\mathcal{\C} \in \mathbb{R}^{l\times m\times p}$;
\STATE 2 $\bar{\A}\leftarrow{\rm{fft}}(\mathcal{A},[\ ],3)$, $\bar{\B}\leftarrow\rm{fft}(\mathcal{B},[\ ],3)$, and $\bar{\C}\leftarrow\rm{fft}(\mathcal{C},[\ ],3)$
\STATE 3 \textbf{for} $i\leftarrow 1 $ \textbf{to} $p$  \textbf{do}
\STATE 4 \ \   $\bar{Y}^{(i)} = \bar{A}^{(i)} \bar{B}^{(i)}$;  $\bar{W}^{(i)} = \bar{C}^{(i)} \bar{A}^{(i)}$;
\STATE 5 \ \  $[\bar{Q}^{(i)}, \bar{R}^{(i)}]=\mbox{qr}(\bar{Y}^{(i)},0)$;
\STATE 6 \ \  $\bar{X}^{(i)}=(\bar{C}^{(i)}\bar{Q}^{(i)})^\dagger \bar{W}^{(i)}$;
\STATE 7 \ \  $\tilde{A}^{(i)}=\bar{Q}^{(i)}\bar{X}^{(i)}$.
\STATE 8 \textbf{end}
\STATE 9 Construct the rank-$k$ approximation $ \hat \A \leftarrow \rm{ifft}(\tilde{\A},[\ ],3)$.
\end{algorithmic}
\end{algorithm}

\begin{algorithm}[H]
\renewcommand{\thealgorithm}{3} 
\caption{Low-Rank Approximation, Fourier domain version}
\label{alg:B}
\begin{algorithmic}
\STATE {\textbf{Input:} $\mathcal{A} \in \mathbb{R}^{m\times n\times p}$; sketch size parameters $k\le l$}
\STATE {\textbf{Output:} low-rank approximation $ \hat \A$}

\STATE 1 Generate Gaussian random tensors $\B \in \R^{n\times k\times p}$, $\mathcal{\C} \in \mathbb{R}^{l\times m\times p}$;
\STATE 2 $\bar{\A}\leftarrow{\rm{fft}}(\mathcal{A},[\ ],3)$, $\bar{\B}\leftarrow\rm{fft}(\mathcal{B},[\ ],3)$, and $\bar{\C}\leftarrow\rm{fft}(\mathcal{C},[\ ],3)$
\STATE 3 \textbf{for} $i\leftarrow 1 $ \textbf{to} $p$  \textbf{do}
\STATE 4 \ \   $\bar{Y}^{(i)} = \bar{A}^{(i)} \bar{B}^{(i)}$;  $\bar{W}^{(i)} = \bar{C}^{(i)} \bar{A}^{(i)}$;
\STATE 5 \ \  $[\bar{Q}^{(i)}, \bar{R}^{(i)}]=\mbox{qr}(\bar{Y}^{(i)},0)$;
\STATE 6 \ \  $[\bar{S}^{(i)}, \bar{T}^{(i)}]=\mbox{qr}(\bar{C}^{(i)}\bar{Q}^{(i)},0)$;
\STATE 7 \ \  $\bar{X}^{(i)}=(\bar{T}^{(i)})^\dagger (\bar{S}^{(i)})^* \bar{W}^{(i)}$;
\STATE 8 \ \  $\tilde{A}^{(i)}=\bar{Q}^{(i)}\bar{X}^{(i)}$.
\STATE 9 \textbf{end}
\STATE 10 Construct the rank$_t$-$k$ approximation $ \hat \A \leftarrow \rm{ifft}(\tilde{\A},[\ ],3)$.
\end{algorithmic}
\end{algorithm}

\section{A Bound For the Frobenius-Norm Error}

In this section we extend Theorem 4.3 of \cite{TYUC17} to tensors.

We first prove a proposition.

\begin{proposition} \label{p5.1}
Let $\G \in \R^{m \times s \times p}$ be a Gaussian random tensor defined by Definition \ref{d2.2}.    Suppose that $\KK \in \R^{s \times n \times p}$, $\B \in \R^{n \times k \times p}$ and $\C \in \R^{l \times \times p}$.    Then
\begin{equation} \label{**}
\EE\|\C * (\G * \KK) * \B\|_F^2 = \|\C\|_F^2\|\B\|_F^2
\end{equation}
Furthermore, if $m > n+1$, then
\begin{equation} \label{***}
\EE\|(\G * \KK)^\dagger\|_F^2 = {n \over m-n-1}.
\end{equation}
\end{proposition}
{\bf Proof}   Let $\F = \C * (\G * \KK) * \B$.   By Lemma \ref{l2.1}, this is equivalent to $\bar F = \bar C \bar G \bar K \bar B$, and
$$\|\C * (\G * \KK) * \B\|_F^2 = {1 \over p} \|\bar C(\bar G\bar K)\bar B\|_F^2.$$
Then $\bar G$ is a standard Gaussian random matrix, and thus, $\bar G\bar K$ is also a standard Gaussian random matrix.   By Proposition A.1 of \cite{HMT11}, we have
$$\EE\|\bar C(\bar G\bar K)\bar B\|_F^2 = \|\bar C\|_F^2 \|\bar B\|_F^2.$$
Then by Lemma \ref{l2.1}, we have (\ref{**}).

Expression (\ref{***}) can be proved similarly.
\qed

Let $\hat \A$ be the rank$_t$-$k$ approximation of $\A$, obtained by Algorithm 3.  We now split the error $\| \A - \hat \A \|_F^2$ to two parts.

\begin{proposition} \label{p5.2}
Let $\A \in \Re^{m \times n \times p}$, and $\hat A$ be the rank$_t$-$k$ approximation of $\A$, obtained by Algorithm 3, with $\X$ and $\Q$ be the intermediate tensors obtained from Algorithm 3, with $\Q^\top * \Q = \I_{kkp}$.   Then
\begin{equation} \label{e5.15}
\|\A - \hat \A \|_F^2 = \|\A - \Q * \Q^\top * \A \|_F^2 + \|\X - \Q^* \A \|_F^2.
\end{equation}
\end{proposition}
{\bf Proof}  As $\hat A = \Q * \X$, we have
\begin{eqnarray*}
&& \|\A - \hat \A \|_F^2 \\
& = & \|\A - \Q * \Q^\top * \A \|_F^2 + \|\Q * \Q^\top * \A - \Q * \X\|_F^2 + \langle \A - \Q * \Q^\top * \A, \Q * \Q^\top * \A - \Q * \X \rangle \\
& = & \|\A - \Q * \Q^\top * \A \|_F^2 + \| \Q * (\Q^\top * \A - \X)\|_F^2 + \langle (\I_{mmp} - \Q * \Q^\top) * \A, \Q * (\Q^\top * \A -  \X) \rangle\\
& = & \|\A - \Q * \Q^\top * \A \|_F^2 + \|\X - \Q^* \A \|_F^2 + {1 \over p}{\rm trace}\left[\bar A^H(I - \bar Q\bar Q^H)\bar Q(\bar Q^H\bar A - \bar X)\right].\\
\end{eqnarray*}
The first part of the last equality is based upon Lemma 3.19 of \cite{KM11}, in the case that $\Q$ is only partially orthogonal, noting that the proof of Lemma 3.19 of \cite{KM11} holds as long as $\Q^\top * \Q = \I_{kkp}$.

Since $\Q^\top * \Q = \I_{kkp}$, we have $\bar Q^H \bar Q = I$.  Then
$$\bar A^H(I - \bar Q\bar Q^H)\bar Q(\bar Q^H\bar A - \bar X) = \bar A^H(\bar Q - \bar Q \bar Q^H \bar Q)(\bar Q^H\bar A - \bar X) = O. $$
We then have (\ref{e5.15}).
\qed

For the first term on the right hand side of (\ref{e5.15}), we have the following theorem.

\begin{theorem} \label{t5.3}
Suppose that $\A \in \R^{m \times n \times p}$.   Let $s$ be a nonnegative integer such that $s < k-1$.  Let $\B \in \R^{n \times k \times p}$ be a Gaussian random tensor defined by Definition \ref{d2.2}, $\Y$ and $\Q$ being calculated by (\ref{e1.2}) and (\ref{e1.3}) respectively.   Then $\Q$ satisfies
$$\EE_\B\|\A - \Q * \Q^\top * \A \|_F^2 \le (1 + f(s, k)) \cdot \tau_{s+1}^2(\A),$$
where $f(s, t) := {s \over t-s-1}$, and $\tau_{s+1}$ is the tail energy defined by (\ref{e3.11}).
\end{theorem}
{\bf Proof}  
Denote $\P = \Q * \Q^\top$. By Lemma \ref{l2.1} and linearity of the expectation, we have
\begin{eqnarray*}
\EE_\B\|\A - \Q * \Q^\top * \A \|_F^2 & = & \EE_\B\|\A - \P * \A \|_F^2\\
 & = & {1 \over p}\left(\sum_{i=1}^p \EE\left\|\bar A^{(i)} - \bar P^{(i)}\bar A^{(i)}\right\|_F^2\right)\\
 & = & {1 \over p}\left(\sum_{i=1}^p \EE\left\|\bar A^{(i)} - \bar Q^{(i)}\left(\bar Q^{(i)}\right)^H \bar A^{(i)}\right\|_F^2\right).
\end{eqnarray*}
By Theorem 10.5 of \cite{HMT11}, we have
$$\EE\left\|\bar A^{(i)} - \bar Q^{(i)}\left(\bar Q^{(i)}\right)^H \bar A^{(i)}\right\|_F^2 \le (1 + f(s, k)) \cdot \tau_{s+1}^2\left(\bar A^{(i)}\right).$$
Thus,
$$\EE_\B\|\A - \Q * \Q^\top * \A \|_F^2 \le {1 \over p} (1 + f(s, k)) \left(\sum_{i=1}^p \tau_{s+1}^2\left(\bar A^{(i)}\right)\right) = (1 + f(s, k)) \cdot \tau_{s+1}^2(\A).$$
The equality here follows from Proposition \ref{p3.4}.
\qed

We now consider the second term on the right hand side of (\ref{e5.15}).    Let $\Q \in \R^{m \times k \times p}$ be the partially orthogonal tensor in Algorithm 3.   Then
$$\left(\overline{\I_{mmp} - \Q * \Q^\top}\right)^{(i)} = I_{m \times m} - \bar Q^{(i)}\left(\bar Q^{(i)}\right)^H,$$
for $i = 1, \cdots, p$.   Hence, there are matrices $\bar P^{(i)} \in \CC^{m \times (n-k)}$ such that
$$\bar P^{(i)}\left(\bar P^{(i)}\right)^H = I_{m \times m} - \bar Q^{(i)}\left(\bar Q^{(i)}\right)^H,$$
for $i = 1, \cdots, p$.   Let $\P \in \R^{m \times (n-k) \times p}$ be the tensor in the real domain by inverse DFT on diag$\left(\bar P^{(1)}, \cdots, \bar P^{(p)}\right)$.  Then by Lemma \ref{l2.1}, we have
$$\P * \P^\top = \I_{mmp} - \Q * \Q^\top.$$

Assume that $\C_1 := \C * \P \in \R^{l \times (n-k) \times p}$ and $\C_2 : = \C * \Q \in \R^{l \times k \times p}$.   Then we have the following proposition.
\begin{proposition} \label{p5.4}
Assume that the tensor tubal rank of $\C_2$ is $k$.  Then
\begin{equation} \label{e5.16}
\X - \Q^\top * \A = \C_2^\dagger * \C_1 * (\P^\top * \A).
\end{equation}
\end{proposition}
{\bf Proof}  By the algorithm design, $\W = \C * \A$.  Since $\P * \P^\top + \Q * \Q^\top = \I_{mmp}$, we have
$$\W = \C * \A = \C * \P * \P^\top * \A + \C * \Q * \Q^\top * \A = \C_1 * (\P^\top * \A) + \C_2 * (\Q^\top * \A).$$
Since the tensor tubal rank of $\C_2$ is $k$, $\C_2^\dagger$ exists.  Right-multiplying the last display in t-product by  $\C_2^\dagger$, we have
$$\C_2^\dagger * \W = \C_2^\dagger * \C_1 * (\P^\top * \A) + \Q^\top * \A.$$
By the algorithm design, we have $\X = \C_2^\dagger * \W$.   The conclusion follows.
\qed

We further have the following proposition.

\begin{proposition} \label{p5.5}
Suppose that $\C \in \R^{l \times m \times p}$ is a Gaussian random tensor, independent from $\B$.  Then
$$\EE_\C[\X-\Q^\top * \A] = \O,$$
and
$$\EE_\C\|\X - \Q^\top * \A \|_F^2 = f(k, l) \cdot \|\A - \Q * \Q^\top * \A \|_F^2,$$
where  $f(s, t) := {s \over t-s-1}$.
\end{proposition}
{\bf Proof}  Observe tensors $\P$ and $\Q$ are partial isometries with orthogonal ranges.   Because of the marginal property of the standard normal distribution, the random tensors $\C_1$ and $\C_2$ are statistically independent standard normal tensors.   The tensor tubal rank of $\C_2 \in \R^{l \times k \times p}$ almost surely is equal to $k$ as $l \ge k$.

Then, taking the expectation of (\ref{e5.16}), we have
$$\EE_\C[\X - \Q^\top * \A] = \EE_{\C_2}\EE_{\C_1}[\C_2^\dagger * \C_1 * \P^\top * \A] = \O.$$
For the first equality, we use the statistical independence of $\C_1$ and $\C_2$ to write the expectation as an iterated expectation.  Then we observe that $\C_1$ is a tensor with zero mean.

Now, taking the expected squared Frobenius norm of (\ref{e5.16}), we have
\begin{eqnarray*}
\EE_\C\|\X - \Q^\top * \A \|_F^2 & = & \EE_{\C_2}\EE_{\C_1}\|\C_2^\dagger * \C_1 *(\P^\top * \A)\|^2_F\\
& = & \EE_{\C_2}[\|\C_2^\dagger\|_F^2 \cdot \|\P^\top * \A\|_F^2]\\
& = & f(k, l)\cdot \|\P^\top * \A \|_F^2.
\end{eqnarray*}
 The two equalities follow from Proposition 5.1.
\qed

Finally, we have the following theorem on the expected error bound of Algorithm 3.

\begin{theorem}
Assume that the sketch parameters satisfy $l > k+1$.   Suppose that random test tensors $\B \in \R^{n \times k \times p}$ and $\C \in \R^{l \times m \times p}$ are drawn independently from the standard normal distribution.   Then the rank$_t$-$k$ approximation $\hat \A$ obtained from (\ref{e1.4}) satisfies
\begin{equation} \label{e5.17}
\EE\|\A - \hat \A\|_F^2 \le (1+f(k, l)) \cdot \min_{\rho < k-1} (1+f(\rho, k)) \cdot \tau_{\rho+1}^2(\A)
={k \over l-k-1} \cdot \min_{\rho < k-1} {k \over k - \rho -1} \cdot \tau_{\rho+1}^2(\A),
\end{equation}
where $\rho$ is a natural number less than $k-1$, $f(s, t) := {s \over t-s-1}$, and the tail energy $\tau_j^2$ is defined by (\ref{e3.11}).
\end{theorem}
{\bf Proof}  Taking the expectation of (\ref{e5.15}), by the linear independence of $\B$ and $\C$, we have
\begin{eqnarray*}
\EE\|\A - \hat \A \|_F^2 & =& \EE_{\B}\|\A - \Q * \Q^\top * \A \|_F^2 + \EE_{\B}\EE_{\C}\|\X - \Q^* \A \|_F^2\\
& = & (1+f(k, l))\EE_{\B}\|\A - \Q * \Q^\top * \A \|_F^2\\
& \le & (1+f(k, l))(1+f(\rho, k)) \cdot \tau_{\rho+1}^2(\A).\\
\end{eqnarray*}
The second equality follows Proposition \ref{p5.5}.  The last inequality follows Theorem \ref{t5.3}.
\qed

\section{Experiments}\label{Experiments}
In this section, numerical experiments are presented to verify the performance of the proposed T-Sketching algorithms (Algorithm 2 and Algorithm 3). We use the following relative error as a metric of the low rank approximation to the input tensor data:
$$
\mbox{Relative Error} :=\frac{\|\A-\hat{\A}_{out}\|^2_F}{\|\A\|^2_F},
$$
where $\A$ and $\hat{\A}_{out}$ are original tensor and estimated low-rank approximation, respectively.
We also employ the peak signal-to-noise ratio (PSNR) defined as
$$
\mbox{PSNR} :== 10 \log_{10}\frac{n_1n_2n_3\|\A\|^2_{\infty}}{\|\A-\hat{\A}_{out}\|^2_F},
$$
and the algorithm running time to evaluate the test methods.

\subsection{Synthetic experiment}
In this subsection, we perform our numerical tests using some synthetic input tensors $\A \in \R^{10^3\times10^3\times10}$ with decaying spectrum.
\begin{enumerate}
  \item {\bf Polynomially decaying spectrum:} These tensors are f-diagonal tensors, their $j_{th}$ frontal slices take the form
  $\A^{(j)}=\mbox{diag}(\underbrace{1,\ldots,1}_{\min(r,j)},2^{-p},3^{-p},4^{-p},\ldots,(n-\min(r,j)+1)^{-p}) \in \R^{n\times n},$
  where $p>0$. In our test, we fix the parameter $r=10$ which controls the rank of the `significant part' of the input tensor, and we consider two cases:
  \begin{itemize}
    \item Slow polynomial decay: $p=1$.
    \item Fast polynomial decay: $p=2$.
  \end{itemize}

  \item {\bf Exponentially decaying spectrum:} The f-diagonal tensors with $j_{th}$ frontal slices taking the form
  $\A^{(j)}=\mbox{diag}(\underbrace{1,\ldots,1}_{\min(r,j)},10^{-q},10^{-2p},10^{-3p},\ldots,10^{-(n-\min(r,j))q}) \in \R^{n\times n}.$
   In our test, we fix the parameter $r=10$ and consider two cases:
  \begin{itemize}
    \item Slow exponential decay: $q=0.25$.
    \item Fast exponential decay: $q=1$.
  \end{itemize}
\end{enumerate}

Figures 1-4 show four examples of how the relative errors, PSNR, CPU time of T-Sketching algorithms and rt\_SVD algorithm (Algorithm 6 in \cite{ZSKA18}), changes with sketch size $k$, respectively. As we can see from the Figures 1-4, the reconstruction results of T-Sketching 1 (Algorithm 2) and T-Sketching 2 (Algorithm 3) are basically the same, except for some slight differences in CPUtime. Compared with rt\_SVD method, accuracy of the approximation results of ``one pass" T-Sketching method is slightly worse, but the running time of T-Sketching method is much shorter. Since rt\_SVD method needs SVD calculation, its running time is 3 to 5 times of T-Sketching method. In particular, as shown in Figures 3-4,  for the input tensor with exponential decaying, the accuracy of T-Sketching method is competitive to the rt\_SVD method. In this case, with less storage and operation, T-sketching method can achieve the similar accuracy of low rank approximation. Additionally, when k reaches a certain value ($k=70$ for the slow exponential decay and $k=30$ for the fast exponential decay), T-Sketching method can obtain the optimal low rank approximation.

\begin{figure}[H]
  \centering
    \includegraphics[scale=0.5]{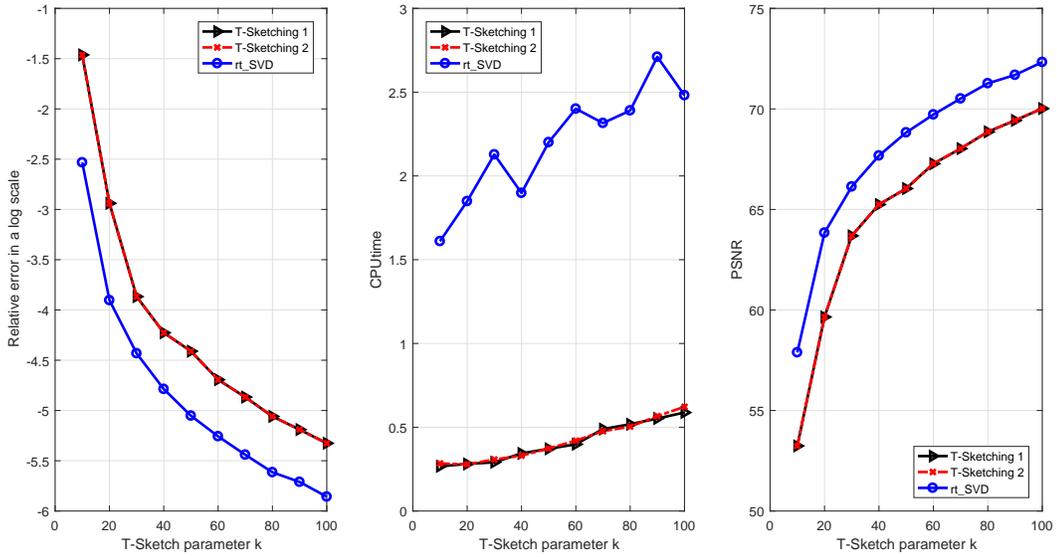}
    \setlength{\abovecaptionskip}{-10pt}
     \setlength{\belowcaptionskip}{-10pt}
  \caption{The Relative Errors, CPUtime, PSNR for the low-rank approximation of synthetic tensor (Slow polynomially decaying spectrum)}
  \label{PolyDecaySlow}
\end{figure}

\begin{figure}[H]
  \centering
    \includegraphics[scale=0.5]{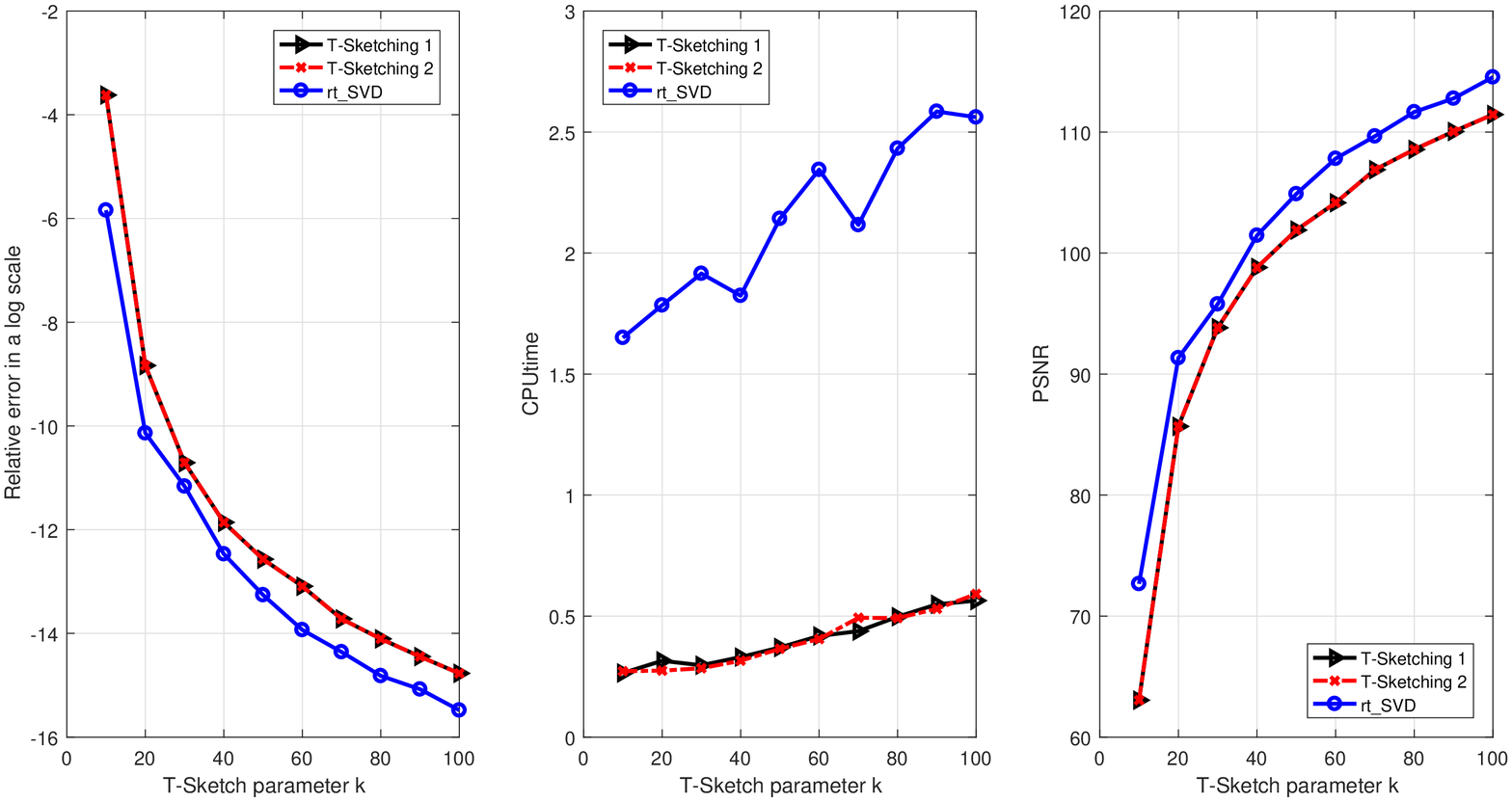}
    \setlength{\abovecaptionskip}{-10pt}
     \setlength{\belowcaptionskip}{-0pt}
  \caption{The Relative Errors, CPUtime, PSNR for the low-rank approximation of synthetic tensor (Fast polynomially decaying spectrum)}
  \label{PolyDecayFast}
\end{figure}

\begin{figure}[H]
  \centering
    \includegraphics[scale=0.5]{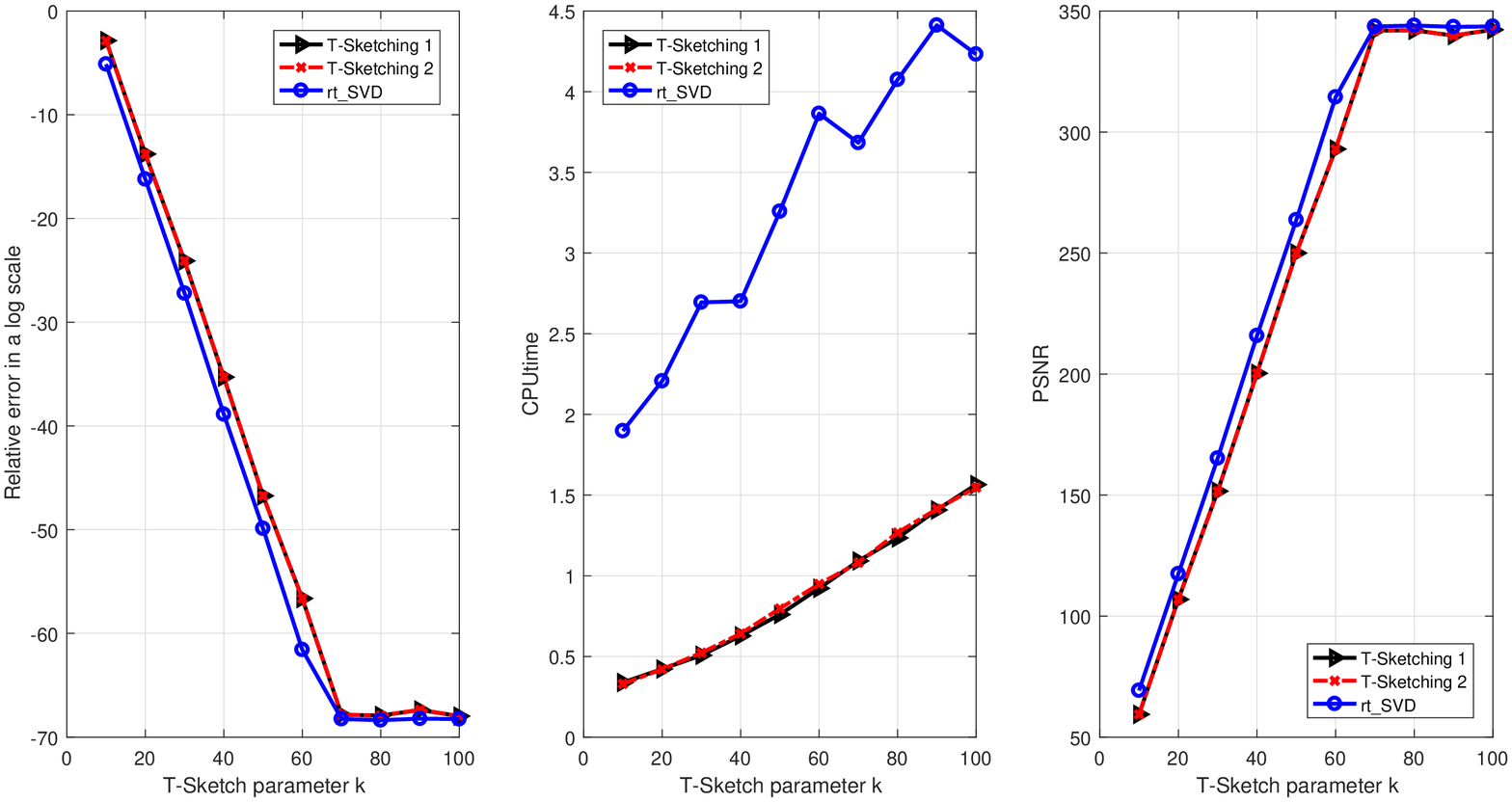}
    \setlength{\abovecaptionskip}{-10pt}
     \setlength{\belowcaptionskip}{-10pt}
  \caption{The Relative Errors, CPUtime, PSNR for the low-rank approximation of synthetic tensor (Slow exponentially decaying spectrum)}
  \label{ExpDecaySlow}
\end{figure}

\begin{figure}[H]
  \centering
    \includegraphics[scale=0.5]{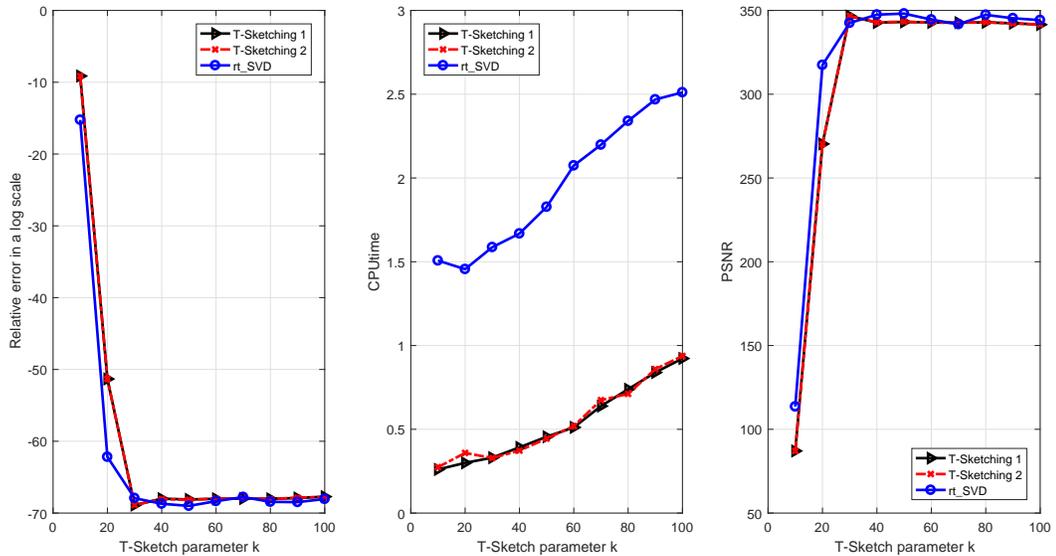}
    \setlength{\abovecaptionskip}{-10pt}
     \setlength{\belowcaptionskip}{-0pt}
  \caption{The Relative Errors, CPUtime, PSNR for the low-rank approximation of synthetic tensor (Fast exponentially decaying spectrum)}
  \label{ExpDecayFast}
\end{figure}

\subsection{T-Sketching on real-world data}

In this subsection, we first test the T-Sketching method for color image data representation and compression via low-rank approximation. The color image, referred to as HDU picture, is of size 1200x1800x3. As shown in Figure 5, the tensor has decaying spectrum. We can see from Figures 6-8 that T-Sketching method can obtain better low rank approximation with the increase of sketch size $k$. In the case of $k = 200$, the relative error of low rank approximation generated by T-Sketching method is less than $10^{-5}$, and the PSNR is 26.27. When $k=600$, the PSNR increased to 38.72 and the the relative error is less than $10^{-8}$. The curve in Figure 8 shows that T-Sketching method and rt-SVD method have similar performance, but T-Sketching method is more effective in terms of running time.

\begin{figure}[H]
  \centering
    \includegraphics[scale=0.5]{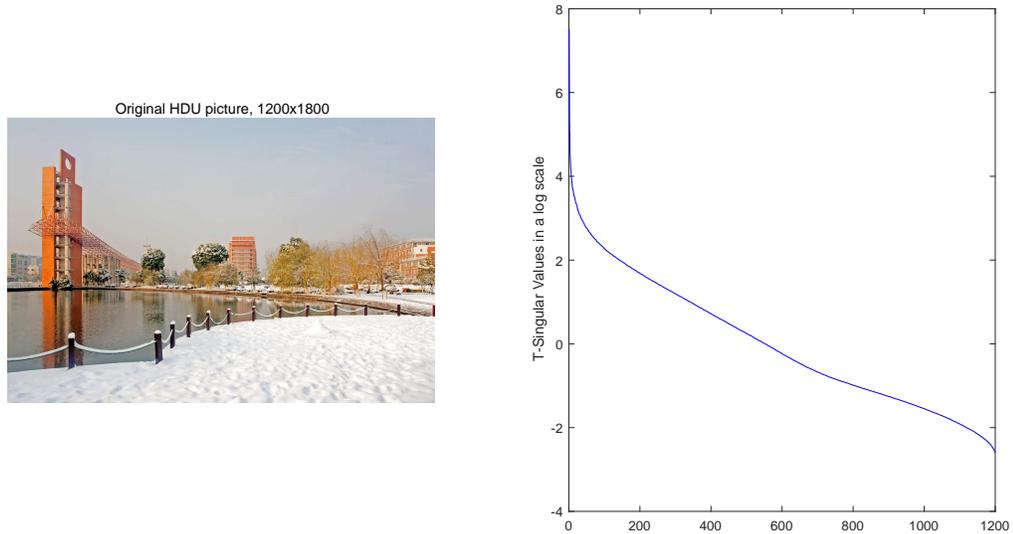}
    \setlength{\abovecaptionskip}{-20pt}
     \setlength{\belowcaptionskip}{-20pt}
  \caption{The original HDU picture and its T-Singular Values}
  \label{HDU-T-Singular-Values}
\end{figure}

\begin{figure}[H]
  \centering
    \includegraphics[scale=0.5]{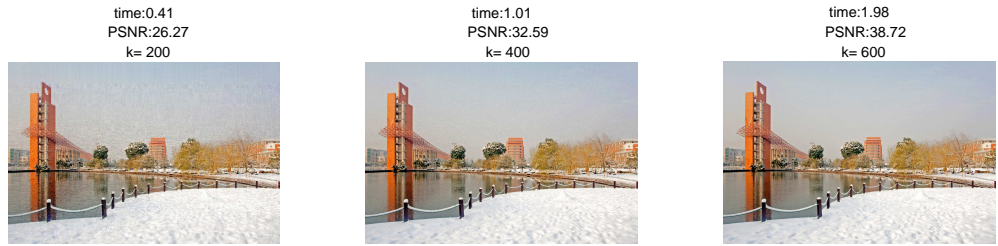}
    \setlength{\abovecaptionskip}{-80pt}
    \setlength{\belowcaptionskip}{0pt}
  \caption{Low-rank approximation via T-Sketching with different Sketch size}
  \label{HDU-T-Singular-Values}
\end{figure}

\begin{figure}[H]
  \centering
    \includegraphics[scale=0.85]{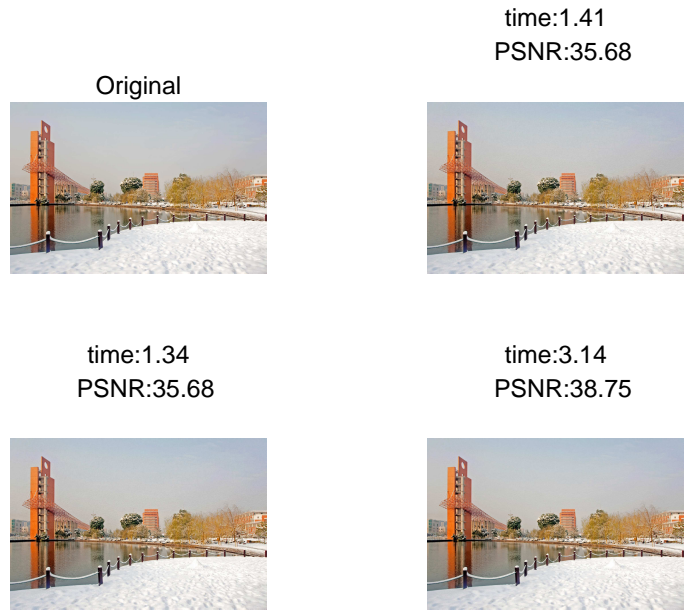}
    \setlength{\abovecaptionskip}{-10pt}
  \caption{Low-rank approximation via T-Sketching v.s. rt-SVD, with sketch size k=500. The upper right and the lower left correspond to the approximation results of T-Sketching 1 and 2, respectively, and the lower right is the approximation results of rt-SVD method.}
  \label{HDU-T-Singular-Values}
\end{figure}

\begin{figure}[H]
  \centering
    \includegraphics[scale=0.5]{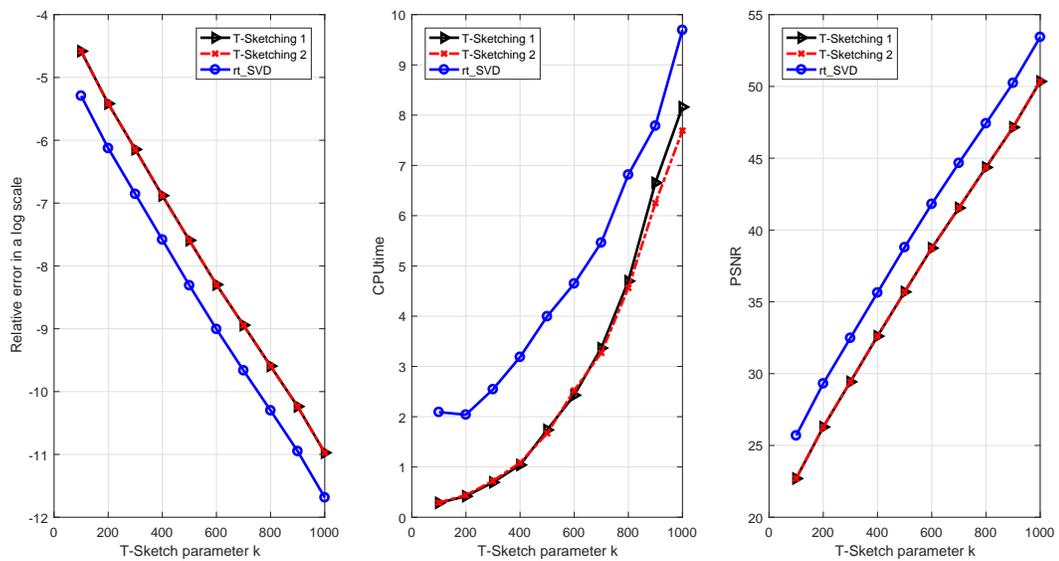}\\
  \caption{The Relative Errors, CPUtime, PSNR for the Low-rank approximation of HDU picture (T-Sketching v.s. rt\_SVD)}
  \label{HDU-T-Singular-Values}
\end{figure}

\begin{figure}[H]
  \centering
    \includegraphics[scale=0.5]{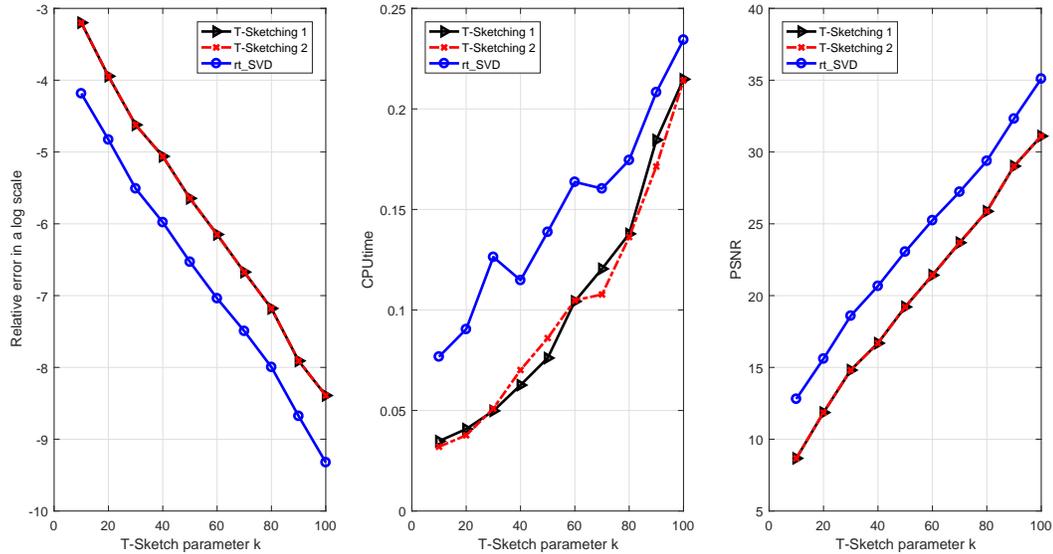}\\
  \caption{The Relative Errors, CPUtime, PSNR for a video (144x176x30 tensor)}
  \label{HDU-T-Singular-Values}
\end{figure}

We also evaluate T-Sketching method on the widely used YUV Video Sequences. Take `hall monitor' video as an example, we only use the first
30 frames. Then the size of the tensor is 144x176x30. As shown in the Figure 9, similar performance can be observed. Therefore, we can conclude that
T-Sketching method avoids computing t-SVD and only needs to store two tensor sketch and two tensors of smaller sizes, leading to a higher algorithm
efficiency.

\bigskip



\bigskip


\begin{thebibliography}{99}

\bibitem{ANW14} H. Avron, H. Nguyen and D. Woodruff, ``Subspace embeddings for the polynomial kernel'', in: Z. Ghahramani, M. Welling, C. Cortes, N.D. Lawrence and K.Q. Weiberger, eds., {\sl Advances in Neural Information Processing Systems \bf 27}, pp. 2258-2266, Curran Associates, Inc., 2014.

\bibitem{Br10} K. Braman, ``Third-Order tensors as linear operators on a space of matrices'', {\sl Linear Algebra and Its Applications \bf 433} (2010) 1241-1253.

\bibitem{CXZ20} Y. Chen, X. Xiao and Y. Zhou, ``Multi-view subspace clustering via simultabeously learing the representation tensor and affinity matrix'', {\sl Pattern Recognition \bf 106} (2020) 107441.

\bibitem{DSSW18} H. Diao, Z. Song, W. Sun and D. Woodruff, ``Sketching for Kronecker product regression and p-splines'', in: A. Storkey and F. Perez-Cruz, eds., {\sl Processings of the Twenty-First International Conference on Artificial Intelligence and Statistics \bf 24}, pp.1299-1308, Canary Island, 2018.

\bibitem{HMT11} N. Halko, P.G. Martinsson and J.A. Tropp, ``Find structure with randomness: Probablistic algorithms for constructing approximate matrix decompositions'', {\sl SIAM Review \bf 53} (2011) 217-288.

\bibitem{KBHH13} M. Kilmer, K. Braman, N. Hao and R. Hoover, ``Third-order tensors as operators on matrices: A theoretical and computational framework with applications in imaging'', {\sl SIAM Journal on Matrix Analysis and Applications \bf 34} (2013) 148-172.

\bibitem{KM11} M. Kilmer and C.D. Martin, ``Factorization strategies for third-order tensors'', {\sl Linear Algebra and Its Applications \bf 435} (2011) 641-658.

\bibitem{KMP08} M. Kilmer, C.D. Martin and L. Perrone, ``A third-order generalization of the matrix svd as a product of third-order tensors'', {\sl Tech. Report \bf TR-2008-4} Tufts University, Computer Science Department, 2008.

\bibitem{LYQX20} C. Ling, G. Yu, L. Qi and Y. Xu, ``A parallelizable optimization method for missing internet traffic tensor data'', arXiv:2005.09838, 2020.


\bibitem{MB18} O.A. Malik and S. Becker, ``Low-rank Tucker decomposition of Large tensors using TensorSketch'', S. Bengio, H. Wallach, H. Larochelle, K. Grauman, N. Cesa-Bianchi and R. Garnett, eds., {\sl Advances in Neural Information Processing Systems \bf 31}, pp. 10117-10127, Curran Associates, Inc., 2018.

\bibitem{MQW20} Y. Miao, L. Qi and Y. Wei, ``Generalized tensor function via the tensor singular value decomposition based on the T-product'', {\sl Linear Algebra and Its Applications \bf 590} (2020) 258-303.

\bibitem{Pa13} R. Pagh, ``Compressed matrix multiplication'', {\sl ACM Transaction on Computation Theory \bf 5} (2013) 1-17.

\bibitem{PP13} N. Pham and R. Pagh, ``Fast and scalable polynomial kernels via explicit future maps'', {\sl Processdings fo the 19th ACM SIGKDD International Conference on Knowledge Discovery and Data Mining \bf KDD13} pp.239-247, 2013.


\bibitem{SHKM14} O. Semerci, N. Hao, M.E. Kilmer and E.L. Miller, ¡°Tensorbased formulation and nuclear norm regularization for multienergy computed tomography¡±, {\sl IEEE Transactions on Image Processing \bf 23} (2014) 1678¨C1693.

\bibitem{SNZ21} G. Song, M.K. Ng and X. Zhang, ``Robust Tensor Completion Using Transformed Tensor SVD'', {\sl Numerical Linear Algebra with Applications} doi.org/10.1002/nla.2299.

\bibitem{TYUC17} J.A. Tropp, A. Yurtsever, M. Udell and V. Cevher, ``Practical sketching algorithms for low-rank matrix approximation'', {\sl SIAM Journal on Matrix Analysis and Applications \bf 38} (2017) 1454-1485.

\bibitem{WTSA15} Y. Wang, H-S.Tung, A.J. Smola and A. Anandkumar, ``Fast and guranteed tensor decomposition via sketching'', in: C. Cortes, N.D. Lawerence, D.D. Lee, M. Sugiyama and R. Garnett, eds',  {\sl Advances in Neural Information Processing Systems \bf 28}, pp. 991-999, Curran Associates, Inc., 2015.

\bibitem{XCGZ21} X. Xiao, Y. Chen, Y.J. Gong and Y. Zhou, ``Low-rank reserving t-linear projection for robust image feature extraction'', {\sl IEEE Transactions on image processing \bf 30}, (2021) 108-120.

\bibitem{XCGZ21a} X. Xiao, Y. Chen, Y.J. Gong and Y. Zhou, ``Prior knowledge regularized multiview self-reprresentation and its applications'', {\sl IEEE Transactions on neural networks and learning systems}, in press.



\bibitem{YHHH16} L. Yang, Z.H. Huang, S. Hu and J. Han, ``An iterative algorithm for third-order tensor multi-rank minimization'', {\sl Computational Optimization and Applications \bf 63} (2016) 169-202.


\bibitem{ZSKA18} J. Zhang, A.K. Saibaba, M.E. Kilmer and S. Aeron, ``A randomized tensor singular value decomposition based on the t-product'',
    {\sl Numerical Linear Algebra with Applications \bf 25} (2018) e2179.

\bibitem{ZA17} Z. Zhang and S. Aeron, ``Exact tensor completion using t-SVD'', {\sl IEEE Transactions on Signal Processing  \bf 65} (2017) 1511-1526.

\bibitem{ZEAHK14} Z. Zhang, G. Ely, S. Aeron, N. Hao and M. Kilmer, ``Novel methods for multilinear data completion and de-noising based on tensor-SVD'', {\sl Proceedings of the IEEE Conference on Computer Vision and Pattern Recognition, ser.  \bf CVPR '14} (2014) 3842-3849.

\bibitem{ZLLZ18} P. Zhou, C. Lu, Z. Lin and C. Zhang, ``Tensor factorization for low-rank tensor completion'',  {\sl IEEE Transactions on Image Processing  \bf 27} (2018) 1152-1163.
































%

%

%

%

%

%

%

%

%

%

%

%

%

%

%

%

%

%

%

%

%

%

%

%

%

\end{thebibliography}
\end{document}